\documentclass[final,hidelinks,onefignum,onetabnum]{siamart251216}

\usepackage{algpseudocode}
\usepackage[sort]{cite}
\usepackage{amssymb}

\usepackage{lipsum}
\usepackage{amsfonts}
\usepackage{graphicx}
\usepackage{epstopdf}
\usepackage{algorithm}
\ifpdf
  \DeclareGraphicsExtensions{.eps,.pdf,.png,.jpg}
\else
  \DeclareGraphicsExtensions{.eps}
\fi

\newsiamremark{remark}{Remark}
\newsiamremark{hypothesis}{Hypothesis}
\crefname{hypothesis}{Hypothesis}{Hypotheses}
\newsiamthm{claim}{Claim}
\newsiamremark{fact}{Fact}
\crefname{fact}{Fact}{Facts}

\headers{}{}

\title{CLUSTER: Derivative-free optimization of smooth functions with parameter-change costs}

\author{
  Serena Landers\textsuperscript{1}
  \and Sahil Pontula\textsuperscript{1,2}
  \and Shiekh Zia Uddin\textsuperscript{3}
  \and Sachin Vaidya\textsuperscript{1,4,5}
  \and Marin Solja{\v{c}}i{\'{c}}\textsuperscript{1,4,5}
  \and Steven G. Johnson\textsuperscript{4,6}
}

\usepackage{amsopn}

\ifpdf
\hypersetup{
  pdftitle={The CLUSTER Algorithm},
  pdfauthor={}
}

\begin{document}
\maketitle

\footnotetext[1]{Department of Physics, Massachusetts Institute of Technology.}
\footnotetext[2]{Department of Electrical Engineering and Computer Science, Massachusetts Institute of Technology.}
\footnotetext[3]{Nokia Bell Labs.}
\footnotetext[4]{Research Laboratory of Electronics, Massachusetts Institute of Technology.}
\footnotetext[5]{NSF Institute for Artificial Intelligence and Fundamental Interactions (IAIFI).}
\footnotetext[6]{Department of Mathematics, Massachusetts Institute of Technology.}

\begin{abstract}

We introduce the CLUSTER algorithm (\textbf{c}oordinate-\textbf{l}evel \textbf{u}pdate \textbf{s}trategy for \textbf{t}rust-region step \textbf{e}valuation \textbf{r}efinement) for local derivative-free optimization problems where there is a cost to changing each parameter (or clusters of parameters). For example, this type of cost model is appropriate for optimizing robot-controlled laboratory experiments, in which a robot may  incur a separate motion for each parameter cluster to be adjusted. We build off of a class of quadratic-interpolation optimization algorithms by Powell and Conn that are known to perform well for twice-differentiable objectives (e.g.~low-noise experiments), and show that the CLUSTER variants improve performance on a variety of test problems (including an optics laboratory experiment) by around 50\%, and greatly outperform common competing algorithms for laboratory optimization (Bayesian optimization and Nelder--Mead). We also adapt the convergence proof of the Conn algorithm to obtain a similar convergence guarantee for CLUSTER-Conn.

\end{abstract}



\section{Introduction}

\begin{figure}
    \centering
    \includegraphics[width=0.5\linewidth]{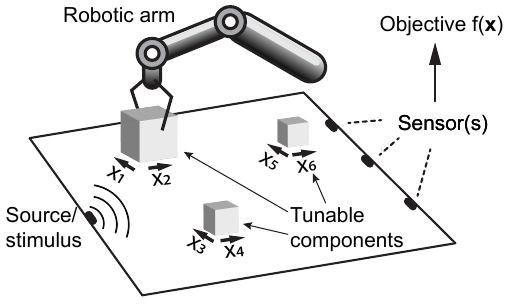}
    \caption{The CLUSTER algorithm is designed for the type of optimization task shown here. The objective is a function of the measurements taken by the sensor(s).}
    \label{fig:fig1}
\end{figure}

The goal of most optimization algorithms is to minimize (or maximize) a function while evaluating the function as few times as possible. For experimental settings, however, where some measured quantity is minimized over a set of physical parameters, the primary cost may not be the function evaluation, but instead the changing of parameters. Consider an example of the generic case shown in Fig.~\ref{fig:fig1}, inspired by Ref.~\citenum{uddin_ai-driven_2025}, where the parameters are the positions of the three optical components that are moved by a robot, and the objective function is the intensity of a laser beam sent through the assembly. In this situation, measuring the light intensity takes considerably less time than moving a component, which requires the robot to physically move to the location of the object and then move it to its new position. Worse, evaluating at an arbitrary point in parameter space requires the movement of all three components and thus takes roughly triple the amount of time of only moving one. This cost model is not limited to optics; in high-throughput autonomous searches over material or chemical compositions~\cite{arnold_cloud_2022, burger_mobile_2020, dai_autonomous_2024, tom_self-driving_2024}, for instance, a robotic platform can often alter only one component at a time but may rapidly measure the resulting material property.

Here, we focus on \emph{local} optimization, in which algorithms can exploit smoothness in an objective function~\cite{ragonneau_model-based_2022}, whereas the number of function evaluations (or parameter moves) required for global optimization may be prohibitive in an experimental setting. Furthermore, in objectives involving experimental measurements, derivative information can be prohibitively expensive or impossible to procure (with rare exceptions~\cite{guillamon_insitu_2025}), necessitating a derivative-free optimization (DFO) algorithm. Whereas previous authors have applied existing DFO algorithms unmodified to experimental parameters~\cite{sun_derivative-free_2020,bezerra_simplex_2016}, one can potentially gain real-world efficiency by modifying the algorithms to target a new metric: minimizing the number of parameter changes.

This work introduces the CLUSTER algorithm (\textbf{c}oordinate-\textbf{l}evel \textbf{u}pdate \textbf{s}trategy for \textbf{t}rust-region step \textbf{e}valuation \textbf{r}efinement) to address this alternative cost metric. Our algorithm can be incorporated as a modification to existing DFO algorithms, and in particular we show its application (CLUSTER-Powell) to the PRIMA implementation~\cite{zhang_2023} of Powell's NEWUOA~\cite{pardalos_newuoa_2006}, a quadratic interpolation trust region algorithm. As NEWUOA lacks a formal convergence proof, we also apply it (CLUSTER-Conn) to a provably convergent variant of this algorithm proposed by Conn~\cite{conn_introduction_2009}. The Powell and Conn algorithms are good starting points: they are known to perform well for local optimization of smooth functions, because they infer both first- and second-derivative information~\cite{powell_uobyqa_2002}. We test both on numerical test functions, and we also demonstrate CLUSTER-Powell on an experimental system to optimize the coupling of a laser beam into a fiber, and find that the CLUSTER version of both hastens convergence by about 50\% compared to the original algorithms. As optimization algorithms are increasingly used in experiments, driven particularly by the advancements in the reliability of robotics performing these optimization tasks~\cite{uddin_collaborative_2025}, algorithms that account for experimentally-relevant cost models can make these processes more efficient.

\subsection{Related Work} A related cost model has been studied in the field of online optimization and multi-armed bandit optimization, where algorithms have been developed that take into account a cost incurred when changing the decision between rounds $t$ and $t+1$, referred to as a switching cost~\cite{dekel_bandits_2014, lin_online_2012, liu_online_2022}, but these algorithms are designed for a particular kind of sequential objective function. Bayesian optimization algorithms have been designed to account for parameters grouped into sequential modules, where changing one module incurs the cost of changing all subsequent modules~\cite{Lin_lambo_2020, torresi_multi-stage_2025}. Similar in spirit to CLUSTER, this metric suits physical optimization problems, as it is targeted to situations in which parameters correspond to steps of an experimental process, so that modifying an early step forces all later steps to be rerun.  

\subsection{Overview of quadratic trust-region algorithms}
\label{ssec:notation}

Quadratic interpolation trust-region algorithms minimize a function $f$ in $n$ dimensions by using a set $X$ of interpolation points to fit a quadratic polynomial $q$, most often with fewer than the $(n+1)(n+2)/2$ points that are required to fully interpolate the polynomial. This class of algorithms performs well for twice-differentiable objectives, corresponding to  noise-insensitive experiments with a suitably-constructed smooth merit function. At each iteration, each algorithm generally follows the steps:

\begin{enumerate}
    \item Fit a quadratic polynomial $q$ to the interpolation set.
    \item Minimize the model within a trust region (a sphere of radius $R$) and evaluate the objective at this minimizer.
    \item Update the interpolation set and trust region radius depending on the agreement between the current model and the true objective function value.
    \item A poor set of points may cause the model to be ill-posed/ill-conditioned, for example if they lie in a lower-dimensional subspace. In such cases, take a model-improvement step: evaluate the objective function at an additional point chosen purely geometrically.
\end{enumerate}

The interpolation model is a quadratic polynomial 
\begin{equation}\label{eq:q} q(x_\text{c}+s)=f(x_\text{c}) + g^Ts+\frac{1}{2}s^T Hs\,,\end{equation} and is denoted $(g, H)$ when it is clear which point $x_\text{c}$ is the center of the model. The center of the model is nearly always the point in $X$ with the minimum function value. At the start of the trust region step, we let $m$ be the number of points in the interpolation set. In Powell's algorithms, $m$ is fixed throughout the algorithm (the default value is $2n+1$), while in Conn's the size of the interpolation set can fluctuate. 

For convergence proofs, one assumes that the function is bounded from below and twice-differentiable with a bounded Hessian~\cite{conn_introduction_2009}. 

For CLUSTER, the parameters $1, \dots, n$ are grouped into disjoint subsets, which we call ``clusters," referring to those parameters that can be moved simultaneously. For example, the two degrees of freedom of one object in Fig.~\ref{fig:fig1} ($x_1$ and $x_2$, for example) can be moved together and therefore fall into the same cluster, while degrees of freedom of two different objects cannot and hence are in different clusters.

\section{Algorithm}
\label{sec:algorithm}
\begin{figure}
    \centering
    \includegraphics[width=\linewidth]{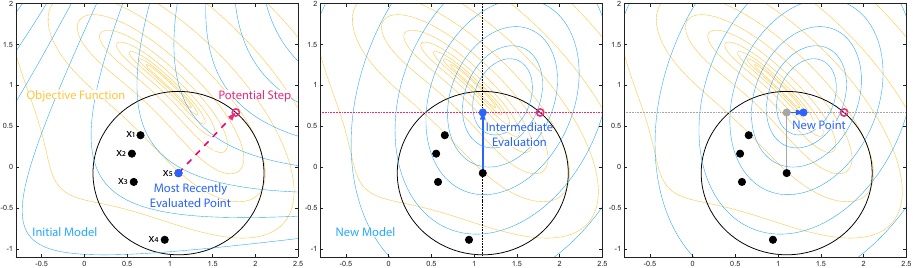}
    \caption{A schematic of Algorithm \ref{alg:evaluate_tr} applied to the function $f = \sum_i \{(x_i-1)^4\} + t^2$, where $t = \mathbf{x}^T A \mathbf{x} - \mathbf{1}^T A \mathbf{1}$ and $A = \begin{pmatrix} 1 & 0.5 \\ 0.5 & 1\end{pmatrix}$. The initial points $X$ are marked as $\mathbf{x}_1, \dots, \mathbf{x}_5$. As is standard in NEWUOA, there are $2n+1 =5$ interpolation points to begin; the additional degree of freedom is fixed by minimizing the Frobenius norm of the Hessian. In each panel, the most recently evaluated point $\mathbf{x}_\text{rec}$ is marked in blue, and the potential point $\mathbf{x}_\text{pot}$ is marked in pink. The contours of the objective function are shown in gold, while the interpolation models, fit to the set $X$ for the first panel and $X$ plus the intermediate evaluation point for the second and third panels, are shown in light blue. }
    \label{fig:fig2}
\end{figure}

The main modification made by CLUSTER lies in how the trust-region step is evaluated, which is described in Algorithm~\ref{alg:evaluate_tr} and depicted in Fig.~\ref{fig:fig2}. When there is only one cluster (i.e., all parameters move together), the original algorithm is recovered. The occasional model-improvement steps are not modified, as the purpose of this type of step is to maximally improve the geometry of the interpolation set, which generally cannot be achieved by changing just some of the clusters. However, as each cluster is changed, it is possible to evaluate at these intermediate points and store the points for future use, or to check for convergence if a minimum target value for the function is set.  There are, in principle, many different ways that these intermediate function evaluations could be exploited, but we have found that Algorithm~\ref{alg:evaluate_tr} works well empirically (as shown in Sec.~\ref{sec:results}) and can also be proved to converge theoretically (Appendix~\ref{asec:convergence}).  In this section, we describe the motivation underlying the choices of Algorithm~\ref{alg:evaluate_tr}.

\begin{algorithm}
    \caption{The trust-region step in the CLUSTER algorithm, which is its main change to the underlying quadratic-model optimization methods.}
    \label{alg:evaluate_tr}
    \begin{algorithmic}[1]

    \State \textbf{Inputs:} Trust region radius $R$; potential point $\mathbf{x}_{\text{pot}}$; most recently evaluated point $\mathbf{x}_{\text{rec}}$; interpolation set $X$ with current center $\mathbf{x}_{\text{c}}$ and associated values of the objective function $f$; model parameters $q_0 = (g_0, H_0)$; hyperparameters $\kappa_1 \ll 1$ and $\kappa_2 \in (0,1)$

    \State \textbf{Initialize:} $q \leftarrow q_0$; $\mathbf{x}_{\text{cur}} \leftarrow \mathbf{x}_{\text{rec}}$; $\mathbf{x}_{\text{best}} \leftarrow \mathbf{x}_{\text{c}}$

    \State Sort coordinate clusters in descending order by the magnitude of $\mathbf{x_\text{pot} - x_{\text{c}}}$ projected onto each cluster.

    \ForAll{cluster $\in \text{sorted clusters}$}

        \State Set $\mathbf{x}_{\text{cur}}[\text{cluster}] \leftarrow \mathbf{x}_{\text{pot}}[\text{cluster}]$
        
        \If{not the first cluster and 
        $\Vert \mathbf{x}_{\text{cur}} - \mathbf{x}_{\text{rec}} \Vert_2 < \kappa_1 R$}
    
            \State $\mathbf{x}_{\text{cur}}[\text{cluster}] = \mathbf{x}_{\text{rec}}[\text{cluster}]$
        
            \State \textbf{continue}
        \EndIf

        \State Evaluate $f(\mathbf{x}_{\text{cur}})$, $\mathbf{x}_{\text{rec}} \leftarrow \mathbf{x}_{\text{cur}}$

        \If{$f(\mathbf{x}_{\text{cur}}) < f(\mathbf{x}_{\text{best}})$}

            \State $\mathbf{x}_{\text{best}} \leftarrow \mathbf{x}_{\text{cur}}$

        \EndIf
            
        \State Set $r = \dfrac{f(\mathbf{x}_{\text{c}}) - f(\mathbf{x}_{\text{cur}})}
          {q(\mathbf{x}_{\text{c}}) - q(\mathbf{x}_{\text{cur}})}$

        \If{$r \in [1-\kappa_2,\, 1+\kappa_2]$ or last cluster}
            \State \textbf{continue}
        \EndIf

        \State Update $q$ by interpolating at $X \cup \{\mathbf{x}_{\text{cur}}\}$, while minimizing the Frobenius norm of the change in $H$.

        \State Recompute $\mathbf{x}_\text{pot}$ by minimizing the model $q$ over the remaining coordinates.

    \EndFor

    \State $\mathbf{x}_{\text{new}} \leftarrow \mathbf{x}_{\text{cur}}$

    \If{$f(\mathbf{x}_{\text{best}}) < f(\mathbf{x}_{\text{c}}) < f(\mathbf{x}_{\text{new}})$}
    
        \State $\mathbf{x}_\text{new} \leftarrow \mathbf{x}_{\text{best}}$

    \EndIf

    \Return{$\mathbf{x_\text{new}}, g, H$}
    \end{algorithmic}
    
\end{algorithm}

During step~2 of the outline given at the end of the introduction, the original algorithms~\cite{pardalos_newuoa_2006, conn_introduction_2009} minimize the model $q$ within the trust region, a hypersphere of radius $R$ about the center of the model. The function is then directly evaluated at the model minimizer, by changing one cluster at a time in some arbitrary order to get from the current configuration in parameter space to the minimizer. In our cost model, however, it does not add significant time to evaluate the function after each cluster is updated. Therefore, instead of blindly changing all of the parameters, CLUSTER changes one cluster at a time, a measurement is taken after each change, and the intermediate information is used to inform changes to subsequent clusters.  

To start, the algorithm takes in the minimizer of $q$ as $\mathbf{x_\text{pot}}$, the potential step (which we may or may not take). This point is updated throughout so that the clusters that have been changed are fixed, while it remains the minimizer of the current model in the space of the remaining clusters. The algorithm also tracks $\mathbf{x_\text{rec}}$, the point where the objective was most recently evaluated, i.e.~the current physical configuration in an experimental system. The remaining inputs describe the current interpolation state: the interpolation set $X$ with its function values, and the model $q$ with center $\mathbf{x}_\text{c}$, which is most often the point in $X$ with the smallest function value.  (In the example of Fig.~\ref{fig:fig2},  $\mathbf{x}_\text{rec} = \mathbf{x}_\text{c}$, but this is not always the case.)

Given these inputs, clusters are first sorted by the magnitude ($\ell^2$ norm) of the initial potential step $\mathbf{x_\text{pot}} - \mathbf{x_\text{c}}$ projected onto the components of each cluster. Other ranking schemes (the projected gradient norm, predicted model decrease, etc.) were tested, but this empirically performs best. The final point $\mathbf{x_\text{new}}$ is then obtained by updating one cluster at a time, in this order, in a loop over the clusters (line~4).

At each iteration of this loop, the components of $\mathbf{x_\text{cur}}$ in the current cluster are set to the corresponding components of the model minimizer $\mathbf{x_\text{pot}}$. If  $\Vert\mathbf{x_\text{cur}} - \mathbf{x_\text{rec}}\Vert_2 < \kappa_1 R$ (we choose a hyperparameter value $\kappa_1 = 10^{-2}$), $\mathbf{x_\text{cur}}$ is reset to $\mathbf{x_\text{rec}}$ and the algorithm moves to the next cluster, except in the first iteration, where it proceeds to the next step regardless of step size. This check (on line~6) avoids negligible parameter changes that are unlikely to affect the function value; the exception for the first cluster (the most important one for this iteration, according to $q$), prevents trust-region steps with no parameter changes.

After each cluster's update, the objective $f$ is evaluated at the new $\mathbf{x_\text{cur}}$ and $\mathbf{x_\text{rec}}$ is updated to $\mathbf{x_\text{cur}}$ (on line~10). If the ratio $r$ of actual to predicted decrease (line~14) lies within $1 \pm \kappa_2$ (we choose a hyperparameter value $\kappa_2 = 0.1$), the algorithm goes directly to the next cluster, since this evaluation is unlikely to improve the model significantly (and can actually degrade its extrapolation ability). Otherwise, the model $q$ is updated (line~18) to interpolate at the original points $X$ and the new point $\mathbf{x_\text{cur}}$, while keeping the center at $\mathbf{x}_\text{c}$. This most often remains an under-determined interpolation problem, so (similar to Powell's original algorithm~\cite{powell_least_2004}) the remaining degrees of freedom are chosen to minimize the Frobenius norm of the change in $H$. Note that intermediate points evaluated during updates earlier in the cluster ranking are not explicitly included in the interpolation set; their effect is instead captured through the minimal-change condition on the Hessian. The computations to update the model are reviewed in Appendix~\ref{asec:update}. 

After the model updates, $\mathbf{x}_\text{pot}$ is updated to be the minimizer of $q$ in the trust region over the space of the remaining coordinates, with the already completed coordinates fixed to their current values.  (This can be algebraically rewritten into an optimization over only the remaining coordinates with a modified trust-region radius, utilizing the same quadratic trust-region solver.)

During the cluster loop, we track $\mathbf{x}_\text{best}$ as the point with the smallest function value so far. At the end of the iteration, if
\[
f(\mathbf{x}_\text{best}) < f(\mathbf{x}_\text{c}) < f(\mathbf{x}_\text{new}),
\]
where $\mathbf{x}_\text{c}$ is the current model center (typically the best initial interpolation point in $X$), we return $\mathbf{x}_\text{best}$ instead of $\mathbf{x}_\text{new}$. \emph{Always} returning $\mathbf{x}_\text{best}$ could harm the interpolation geometry because it often shares several coordinates with a point in $X$, leading to some coordinates with few distinct values, which requires many model-improvement steps to repair. Thus, $\mathbf{x}_\text{best}$ is returned only when the above condition holds, in which the objective decrease justifies the potential geometry degradation.

The center $\mathbf{x}_\text{c}$ is unchanged during the cluster loop even if an intermediate point results in a smaller objective function value, to ensure the trust region within which we minimize the function does not change throughout the process. Upon entering the algorithm, the trust region is the hypersphere of radius $R$ around the point $\mathbf{x}_\text{c}$, which should remain fixed.  (After the trust-region step, $\mathbf{x}_\text{c}$ is updated by step~3 of the algorithm in Sec.~\ref{ssec:notation}.) 

\section{Results}
\label{sec:results}

This section discusses results for the CLUSTER-Powell algorithm (and competing algorithms) on several numerical test problems, as well as an experimental problem where function values are obtained from optical measurements. The results for CLUSTER-Conn (which does slightly worse in practice for our tests) and its convergence analysis, as well as a discussion of the algorithmic differences, are given in Appendix \ref{asec:conn_cluster}.   Like other work on derivative-free algorithms, we consider optimization problems in moderate dimensions ($\le 30$); for high-dimensional numerical optimization one normally should exploit gradients, whereas in an experimental setting the time requirements for high-dimensional derivative-free optimization is likely to be prohibitive.

\subsection{Numerical test problems}

\begin{figure}
    \centering
    \includegraphics[width=\linewidth]{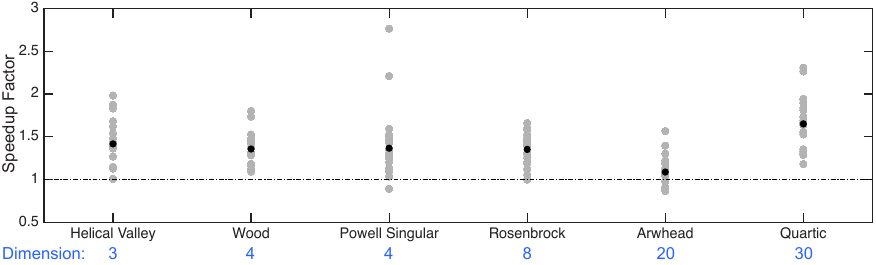}
    \caption{The CLUSTER-Powell algorithm's speed-up factor (in parameter changes to convergence) compared to (non-CLUSTER) NEWUOA~\cite{pardalos_newuoa_2006} (as implemented in PRIMA~\cite{zhang_2023}) for twenty random starting points on six test functions. The geometric means are marked with black points: all $> 1$, and typically 30--50\% speedup.}
    \label{fig:fig3}
\end{figure}

For each test problem, we use 20 random starting points drawn uniformly from a box of side length $2R_0$ centered at the standard starting point, where $R_0$ is also the initial trust-region radius. For each starting point, we run the original (non-CLUSTER) algorithm and also CLUSTER-Powell with each coordinate in its own cluster. We compare them via the speed-up factor, defined as the ratio of the number of parameter changes to those of the original algorithm needed to converge within a given tolerance (defined in Appendix~\ref{asec:test_functions}). A parameter change is whenever any one parameter is modified, regardless of whether the function is evaluated immediately; thus, for the original non-CLUSTER algorithm on an $n$-dimensional problem where all of the parameters are changed at once, a single function evaluation requires $n$ parameter changes. The test functions, mostly taken from the optimization literature, are defined in Appendix \ref{asec:test_functions}, and results are presented in Fig.~\ref{fig:fig3}. CLUSTER tends to provide the greatest benefit in objectives where variables are coupled (mixed derivatives are nonzero), since intermediate evaluations are more informative, but even for nearly uncoupled (Arwhead) and fully uncoupled (quartic) functions, its average speed-up factor exceeds~1.

For illustration purposes, we also compare NEWUOA and CLUSTER-Powell to other popular derivative-free optimization algorithms, Nelder--Mead~\cite{NLopt, NELDERMEAD} and Bayesian optimization~(BO)~\cite{bayesian_opt}, although it is widely known that NEWUOA-style algorithms tend to be superior for \emph{local} optimization of \emph{smooth} functions to methods like Nelder--Mead, BO, or genetic algorithms that do not infer derivative information~\cite{conn_introduction_2009}. The convergence histories for the Helical Valley and Arwhead functions (defined in Appendix \ref{asec:test_functions}) are given in Fig.~\ref{fig:fig4}. BO finds the global optimum by balancing exploration of the space and exploitation of current minima. Due to this procedure, BO's convergence history exhibits wild oscillations, and therefore we show its best-so-far point at each parameter change. As BO is a global algorithm, a finite region over which to search must be specified; we choose a box with side length $4 R_0$, centered at the optimum, where $R_0$ is the initial trust-region radius supplied to CLUSTER-Powell and NEWUOA. Because it is a global-search algorithm, BO is not strictly comparable to local optimization methods like CLUSTER, but we include it because of its popularity in experimental settings in order to quantify the price one pays for its globality. As expected, Nelder--Mead and BO both converge far more slowly than both NEWUOA and CLUSTER-Powell, especially in higher dimensions. In the Supplementary Information, we also compared CLUSTER-Powell and NEWUOA with several coordinate-search methods~\cite[p.~115]{conn_introduction_2009} and found that, particularly for complex functions and in cases where the starting point is far from the optimum, coordinate-search is over an order of magnitude slower for simple coordinate-search algorithms~\cite{conn_introduction_2009} and more than a factor of two slower for the NOMAD hybrid direct-search algorithm~\cite{audet_algorithm_2022}.

\begin{figure}
    \centering
    \includegraphics[width=\linewidth]{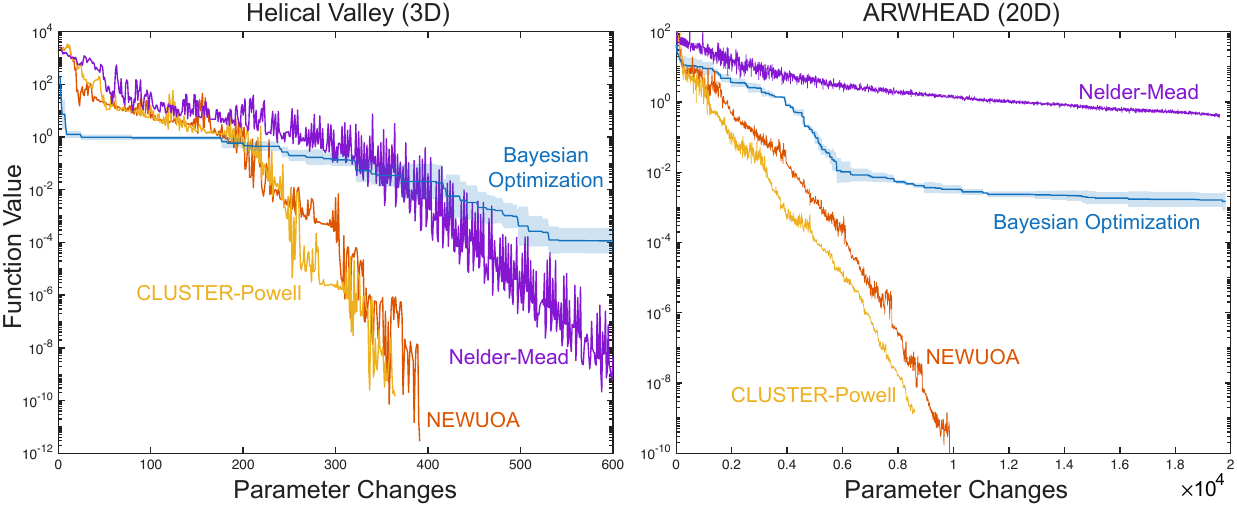}
    \caption{Comparison between CLUSTER-Powell, NEWUOA, Nelder--Mead, and Bayesian optimization for the helical valley and ARWHEAD functions. The Bayesian optimization curve shows the best-so-far value, and is the average of three runs with different sets of randomly-selected initialization points.}
    \label{fig:fig4}
\end{figure}

\subsection{Optimization of a laboratory experiment}
\label{sec:experiment}

\begin{figure}
    \centering
    \includegraphics[width=0.6\linewidth]{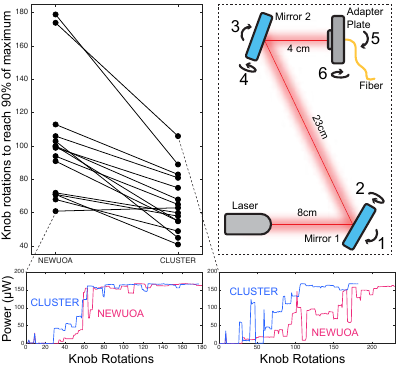}
    \caption{Testing of the CLUSTER-Powell algorithm on optimization of a laboratory experiment (an optical laser--fiber coupler). \emph{Upper~right:} schematic of the experiment.  \emph{Upper~left:} results of 15 trials with different random starting points. \emph{Bottom:} convergence history of two representative trials.  The geometric mean of the speedup factors for CLUSTER (compared to the original NEWUOA algorithm) is $1.5\times$.}
    \label{fig:fig5}
\end{figure}

CLUSTER-Powell was also tested on a physical laboratory experiment for which the objective was to maximize the laser power coupled into a fiber from free space with respect to angular parameters of optical components, inspired by Ref.~\citenum{uddin_ai-driven_2025} and Fig.~\ref{fig:fig1}. A schematic of the configuration is shown in Fig.~\ref{fig:fig5}; the two mirrors and the coupler (adapter plate) each have two knobs controlling tip/tilt, resulting in six total parameters. Fig.~\ref{fig:fig5} presents the results of fifteen optimization runs, where the starting points were randomly selected from a uniform distribution of $\pm 0.12^{\circ}$ centered around the optimum (roughly the maximum angular deviation of the first mirror before the coupled power becomes negligible). Each starting point was tested twice, first with the original NEWUOA algorithm (in which every knob is changed on every function evaluation) and again with CLUSTER-Powell (for one knob per cluster). The left panel shows the number of knob turns (parameter changes) required to reach 90\% of the maximum power eventually achieved (by either algorithm) for that starting point. Two representative convergence histories are given underneath, with dashed lines indicating to which trials they correspond. Across these fifteen trials, the speedup factor (reduction in parameter changes) due to CLUSTER was approximately $1.5 \pm 0.3$ (as computed by the geometric mean, or equivalently the mean and variance on a log scale). We also attempted Bayesian optimization and Nelder--Mead in this experiment. BO exhibited a huge variation in the number of iterations required depending on whether a random sample happened to lie near the optimum, but was typically at least $2\times$ slower. Nelder--Mead was also approximately $2 \times$ slower than NEWUOA, with larger slowdowns when initial starting points were further from the optimum. (Plots of experimental data for BO and Nelder--Mead can be found in the Supplementary Information.)

\section{Conclusion}
We introduced the CLUSTER modification to existing derivative-free optimization (DFO) quadratic trust-region algorithms that minimizes parameter changes rather than function evaluations. Tests show that the algorithm, by using intermediate function evaluations to refine quadratic polynomial models, typically reduces the number of parameter changes needed for convergence by 20--60\%, for both numerical test functions and a physical optimization system. 

There are many opportunities for future work in adapting existing optimization algorithms (and constructing new ones) for experimental settings where the cost model is very different from that of conventional optimization theory. Here, we considered problems where the cost is dominated by parameter changes (perhaps grouped in ``clusters''), where each parameter change has an equal cost, but the algorithm could be extended to parameters with different switching costs, possibly by sorting parameters in descending order of cost rather than by step length (line~3) and by using a stricter evaluation cutoff for parameters with higher switching costs (line~6).If measurements can be made nearly \emph{continuously} as a parameter is adjusted, there is also the possibility of some form of exact line search, either with classic cyclic coordinate descent~\cite{bertsekas_nonlinear_1999} or using the quadratic model to choose search directions.

There have also been efforts on problems where the objective function is broken into ``experimental stages/modules'' that have nonuniform costs (and perhaps yield additional intermediate data that could be exploited). Work on the latter class of problems has mainly focused on modifying Bayesian optimization (BO)\cite{Lin_lambo_2020, torresi_multi-stage_2025}. If the objective function is smooth (not including noise), and a local optimum suffices, then DFO algorithms like Powell's that exploit smoothness can be far more efficient than BO, and it is potentially attractive to extend algorithms like CLUSTER to other cost models such as ``staged'' experiments. Conversely, BO algorithms may be more suitable for problems with discrete parameters and many poor local minima, but they have not yet been adapted to the parameter-switching cost model. 

In experimental settings with significant noise, smoothness in the underlying expected value can still be exploited by extending ideas from stochastic DFO~\cite{larson_derivative-free_2019, deng_adaptation_2006} (analogous to stochastic gradient descent for noisy optimization in machine learning), and it would be interesting to incorporate such ideas into CLUSTER. However, we believe that CLUSTER represents an important first step in extending smooth local DFO to a laboratory-friendly cost model.

\begin{appendix}

\section{Updating the temporary model}
\label{asec:update}

During each iteration of Algorithm \ref{alg:evaluate_tr} when $\mathbf{x}_\text{new}$ is sufficiently different from $\mathbf{x}_\text{rec}$ and the ratio between the predicted decrease and the actual decrease is not close to 1, the model $q$, is updated to interpolate at $X \cup \{\mathbf{x}_\text{new}\}$, with the degrees of freedom fixed by a minimal change condition on the Hessian. 

The update is obtained by solving the quadratic programming (QP) minimization problem:
\[
\min_{\Delta c, \Delta g, \Delta H} \|\Delta H\|_{\text{F}}^2 \quad \text{s.t.} \quad q_{\text{new}}(\mathbf{x}_i) = f(\mathbf{x}_i) \ \forall\, \mathbf{x}_i \in X \cup \{\mathbf{x}_\text{new}\},
\]
where
\[
q_\text{new}(\mathbf{x}) = q_\text{old}(\mathbf{x}) 
+ \Delta c + \Delta g^T(\mathbf{x}-\mathbf{x_\text{c}}) 
+ \tfrac{1}{2} (\mathbf{x}-\mathbf{x_\text{c}})^T \Delta H(\mathbf{x}-\mathbf{x_\text{c}}),
\]
and
\[
q_\text{old}(\mathbf{x}) = c + g^T(\mathbf{x}-\mathbf{x_\text{c}}) + \tfrac{1}{2} (\mathbf{x}-\mathbf{x_\text{c}})^T H (\mathbf{x}-\mathbf{x_\text{c}}).
\]
Although many existing software packages and algorithms can easily solve a small QP like this, the solution here can be expressed compactly in terms of the Lagrange dual problem, similar to the method described by Powell~\cite{powell_least_2004}. The constant $c$ is the value of the function at the center of the model; because the model center is unchanged, $\Delta c = 0$. $q_\text{old}$ already interpolates the objective at the points $X$, so the minimization problem is equivalent to
\[\min_{\Delta c, \Delta g, \Delta H} \|\Delta H\|_{\text{F}}^2 \quad \text{subject to}\]
\[
\begin{aligned} 
\Delta c + \Delta g^T (\mathbf{x}_i - \mathbf{x}_{\text{c}}) &+ \frac{1}{2} (\mathbf{x}_i - \mathbf{x}_{\text{c}})^T \Delta H (\mathbf{x}_i - \mathbf{x}_{\text{c}})
\\ &= \begin{cases} 0, & \mathbf{x}_i \in X, \\[4pt] f(\mathbf{x}_{\text{new}}) - q_{\text{old}}(\mathbf{x}_{\text{new}}),& \mathbf{x}_i = \mathbf{x}_{\text{new}} \end{cases}
\end{aligned}
\]
The Lagrangian of this problem is 
\begin{multline*} 
\mathcal{L} = \frac{1}{2} \sum_{i = 1}^n \sum_{j = 1}^n [\Delta H]_{i,j}^2 \\ - 
\sum_{k = 1}^{m+1} \lambda_k \bigg\{ \Delta c + \Delta g^T(\mathbf{x}_k - \mathbf{x}_{\text{c}}) + \frac{1}{2}(\mathbf{x}_k - \mathbf{x}_{\text{c}})^T \Delta H (\mathbf{x}_k - \mathbf{x}_{\text{c}}) \\ - \big(f(\mathbf{x}_\text{new}) - q_\text{old}(\mathbf{x}_\text{new})\big)\delta_{k, m+1} \bigg\},
\end{multline*}
where the $\lambda_k$'s are Lagrange multipliers. Differentiating with respect to $\Delta c, \Delta g$, and $\Delta H$ results in the KKT equations
\[
\sum_{k=1}^{m+1} \lambda_k (\mathbf{x}_k - \mathbf{x}_\text{c}) =0,  \quad \sum_{k=1}^{m+1} \lambda_k = 0, \] 
\[\Delta H = \frac{1}{2} \sum_{k=1}^{m+1} \lambda_k (\mathbf{x}_k - \mathbf{x}_\text{c})(\mathbf{x}_k - \mathbf{x}_\text{c})^T  = \frac{1}{2} \sum_{k=1}^{m+1} \lambda_k \mathbf{x}_k \mathbf{x}_k^T
\]
We formulate the following matrix equation to solve for the $\lambda_k$'s, $\Delta c$, and $\Delta g$ 
\[ \left(\begin{array}{c|c} A & Y^T \\ \hline Y & 0
\end{array} \right) \begin{pmatrix} \vec{\lambda} \\ \Delta c \\ \Delta g\end{pmatrix} = \big(f(\mathbf{x}_\text{new}) - q_\text{old}(\mathbf{x}_\text{new})\big)\mathbf{e}_{m+1} \, ,
\]
where $\mathbf{e}_{m+1} \in \mathbb{R}^{n+m+2}$  is the $(m+1)$-th Cartesian unit vector, $[A]_{ij} = \frac{1}{4} \{(\mathbf{x}_i-\mathbf{x}_\text{c})^T(\mathbf{x}_j-\mathbf{x}_\text{c})\}^2$ for $i, j = 1, \dots, m+1$ and \[Y = \begin{pmatrix} 1 & 1 & \dots & 1 & 1 \\ 
\mathbf{x}_1 & \mathbf{x}_2 & \dots & \mathbf{x}_{m} & \mathbf{x}_\text{new} \end{pmatrix}\]
Solving this system provides $\Delta c$ and $\Delta g$, and the $\lambda_k$'s give $\Delta H$, yielding the model update: $c \leftarrow c + \Delta c, \quad g \leftarrow g + \Delta g, \quad H \leftarrow H + \Delta H$.

\section{Test functions}
\label{asec:test_functions}
The following test functions were used in the numerical experiments discussed in Section \ref{sec:results}. For each test function, the starting point $\mathbf{x}_0$, the initial trust-region radius $R_0$, and the convergence threshold for computing the speedup factor are given. The twenty trials start at a point randomly selected from the uniform distribution in the $n$-dimensional box $\mathbf{x}_0 \pm R_0$. The same test functions are used for the Conn version of the algorithm introduced in the next section. 

\begin{itemize}
\item \textit{Helical Valley (3D)~\cite{more_testing_1981}}
\begin{align*}f = 100(x_3 - 10\Theta(x_1,x_2))^2 + 100((x_1^2+x_2^2)^{1/2}-1)^2 + x_3^2, \\ \text{where }\Theta(x_1, x_2) = \begin{cases} \tan^{-1}(x_2/x_1)/(2\pi) & x_1 > 0 \\  \tan^{-1}(x_2/x_1)/(2\pi) + 0.5  & x_1 < 0 \\ 0.25 & x_1 =0, x_2 > 0\\
-0.25 & x_1 = 0, x_2 < 0\end{cases}\end{align*}

$\mathbf{x}_0 = (-1, 0, 0), R_0 = 0.5$, Threshold = $10^{-7}$
    
\item \textit{Wood (4D)~\cite{more_testing_1981}:}
\begin{align*} f = \sum_{n=1}^6 f_n, \quad
    f_1 = 100 (x_2 - x_1^2)^2,
    f_2 = (1 - x_1)^2,
    f_3 = 90 (x_4 - x_3^2)^2, \\
    f_4 = (1 - x_3)^2,
    f_5 = 10(x_2 + x_4 - 2)^2, 
    f_6 = 0.1(x_2 - x_4)^2 
\end{align*}

$\mathbf{x}_0 = (-3, -1, -3, -1), R_0 = 1.0$, Threshold = $10^{-7}$

\item \textit{Powell singular (4D)~\cite{more_testing_1981}:}
\begin{align*} 
    f = (x_1 + 10x_2)^2 + 5(x_3 - x_4)^2 + (x_2 -2x_3)^4 + 10(x_1-x_4)^2
\end{align*}

$\mathbf{x}_0 = (-3, 1, 0, 1), R_0 = 1.0$, Threshold = $10^{-6}$

\item \textit{Extended Rosenbrock (8D)~\cite{more_testing_1981}:} $\sum_{n=1}^4 100(x_{2n} - x_{2n-1}^2)^2 + (1-x_{2n-1})^2$

$\mathbf{x}_0 = (-1.2, 1.0, -1.2, 1.0, -1.2, 1.0, -1.2, 1.0), R = 0.5$, Threshold = $10^{-7}$

\item \textit{Arwhead (20D)~\cite{pardalos_newuoa_2006}:}
    $\sum_{n=1}^{19} \Big\{ (x_n^2+x_{20}^2)^2 - 4x_n+3\Big\}$

$\mathbf{x}_0 = (1, 1, \dots, 1, 1, 0), R_0 = 0.5$, Threshold = $10^{-7}$

\item \textit{Quartic (30D):}$ \sum_{n = 1}^{30} (x_n-1)^4$

$\mathbf{x}_0 = (0, 0, \dots, 0, 0), R_0 = 0.5$, Threshold = $10^{-15}$
\end{itemize}

\section{CLUSTER-Conn}
\label{asec:conn_cluster}

This section presents the results for the CLUSTER-Conn algorithm, based on Conn's variant of the quadratic DFO algorithm~\cite{conn_global_2009}. Conn's version of the algorithm has several differences from NEWUOA, which we incorporate into CLUSTER-Conn. Firstly, for convergence-proof purposes, there is a criticality step at the start of every iteration~\cite{conn_introduction_2009} (prior to the trust-region step of Algorithm~\ref{alg:evaluate_tr}). The purpose of Conn's criticality step is to ensure that when the gradient norm (the measure of first-order stationarity) is below some tolerance $\epsilon_{\text{c}}$, it is due to convergence to a local minimum and not a poor model. One iteration of the criticality step consists of improving the model (adding or replacing points as in the model-improvement algorithm) until it is fully linear within the trust region, and shrinking the trust region by a constant factor. This is repeated until the trust-region radius is smaller than some proportion of the gradient norm, $\mu \Vert g\Vert $, which can occur with or without the gradient norm increasing. The former indicates that the model was repaired to the point that the gradient is no longer first-order critical. After exiting the criticality step, it is possible that the convergence criteria for the trust-region radius is met. Otherwise, the algorithm continues and another trust-region step is computed. The criticality step thus ensures the model is fully linear at this stage, and the model's gradient norm is a good measure of stationarity. 

We choose the (dimensionful) parameters $\epsilon_{\text{c}}$ and $\mu$ as follows, since values are not suggested by Conn. The criticality step---which is only useful near convergence of the optimizer---should be triggered when the model's gradient norm $\Vert g\Vert$ is slightly larger than the user’s convergence threshold $g_{\text{tol}}$. That is,
\[
\epsilon_{\text{c}} = r_1 \, g_{\text{tol}},
\]
with $r_1>1$ (e.g., $r_1 = 10$).  When the criticality steps terminate (for $R < \mu \Vert g\Vert$), one would like the trust-region radius $R$ to be close to, but somewhat larger than, the convergence threshold $R_{\text{tol}}$ for $R$; that is,
\[
R \approx r_2 R_{\text{tol}},
\]
for another order-unity constant $r_2>1$. So, when the criticality steps exits we should have $R \approx \mu \Vert g\Vert \approx r_2 R_{\text{tol}}$, while the model-improvement rarely changes the gradient much from its value of $\Vert g\Vert \approx r_1 g_{\text{tol}}$, and combining these approximate relations suggests the hyperparameter choice
\[
\mu \approx \frac{r_2}{r_1} \,\frac{R_{\text{tol}}}{g_{\text{tol}}} \approx \frac{R_{\text{tol}}}{g_{\text{tol}}}.
\]

Another difference between Conn's algorithm and NEWUOA is the construction of the model after adding a new point (i.e., after the trust-region step). Unlike NEWUOA, which updates the model by minimally changing the Hessian while interpolating at the new point, Conn at each step rebuilds a new minimum-Frobenius-norm model. This discards some information from previous iterations, but is compensated for by scanning a ``bank'' of all previously-evaluated points within a radius $tR$ (with $t \geq 1$) from the model center and adding points to the interpolation set $X$ until either $\frac{(n+1)(n+2)}{2}$ points are collected (enough to fully interpolate a quadratic) or the next point would be geometrically unfavorable. In CLUSTER-Conn, the intermediate points (from individual parameter-cluster changes in Algorithm~\ref{alg:evaluate_tr}) are also included in this bank of previously-evaluated points for subsequent model updates, after the trust-region step. Also, the model update in the cluster loop of Algorithm \ref{alg:evaluate_tr} (line~18) is skipped in the (rare) event that $X$ contains $\frac{(n+1)(n+2)}{2}$ points (so that $q$ is already fully determined). 

Our CLUSTER-Conn algorithm incorporates these differences, and also adds new checks to the CLUSTER trust-region step (Algorithm~\ref{alg:evaluate_tr}) to facilitate our convergence proof, as described in the next Appendix. 

Although NEWUOA does not have a global convergence proof, it is more widely used than the Conn's algorithm, probably because of the difference in software availability (we were forced to re-implement the Conn algorithm ourselves), but we also find that NEWUOA tends to be faster empirically in many (but not all) cases. To compare the two implementations of the algorithm, we compute the average cost (in parameter changes) of the original version of the Conn algorithm and CLUSTER-Conn, both relative to NEWUOA, for the same test problems in Appendix~\ref{asec:test_functions}. We also include the cost of CLUSTER-Powell relative to NEWUOA (the inverse of the values in Fig.~\ref{fig:fig3}). 

\begin{table}[h]
    \centering
    \begin{tabular}{|c|c|c|c|c|c|c|}
    \hline
        Cost (relative & Helical & Wood & Powell & Rosenbrock & Arwhead & Quartic \\
         to NEWUOA) & Valley & & Singular & & & \\
        \hline
        Conn & 2.33 & 2.86 & 1.78 & 4.52 & 2.04 & 0.66\\
        \hline
        CLUSTER-Conn& 1.57 & 1.74 & 0.61 & 1.10 & 0.51 & 0.18 \\
        \hline
        CLUSTER-Powell& 0.71 & 0.74 & 0.73 & 0.74 & 0.92 & 0.61\\
        \hline
    \end{tabular}
    \caption{The geometric-mean cost (number of parameter changes) of the Conn algorithm, CLUSTER-Conn, and CLUSTER-Powell, all relative to NEWUOA (the last row is the inverse of the values plotted in Fig.~\ref{fig:fig3}). A relative cost greater than~1 indicates NEWUOA is faster.}
    \label{tab:newuoa_conn_ratios}
\end{table}

CLUSTER-Conn has two advantages over the original Conn algorithm: similar to CLUSTER-Powell it exploits intermediate single-cluster steps for model updates during the trust-region step, but it \emph{also} includes these intermediate evaluations to the bank of points available for the construction of the model during the next iteration.  To separate these two effects, we also considered a variant of Conn that incorporates intermediate evaluations into the bank, but does \emph{not} use them for model updates during the trust-region step.   We found that this indeed produced some marginal speedups of the Conn algorithm (by up to $\approx 30$\% in the higher-dimensional test problems), but was still significantly slower than CLUSTER-Conn.

\section{Convergence analysis}
\label{asec:convergence}

The proof showing convergence of the Conn algorithm to a first-order critical point presented in Refs.~\citenum{conn_introduction_2009, conn_global_2009} can be applied to CLUSTER-Conn with minor modifications described in this appendix. To preserve the convergence guarantee, we introduce several checks beyond the procedure in Algorithm \ref{alg:evaluate_tr}, with fallbacks to the original (full-coordinate) algorithm. Empirically, these fallbacks are rare and do not seem to add a significant number of additional parameter changes. 

The three proof components that must be addressed are: (1) maintaining a well-poised interpolation set, (2) correctly triggering the criticality step and preserving its guarantees, and (3) satisfying the sufficient decrease condition in the trust-region step.

\subsection{Well-poisedness}
Let $q$ be an interpolation for the function $f$. Full linearity~\cite[p.~90]{conn_introduction_2009} in the trust region $\mathbf{x}_k + \mathbf{s}$ for $\Vert \mathbf{s} \Vert \le R$ requires that 
\[
\|\nabla f(\mathbf{x}_k + \mathbf{s}) - \nabla q(\mathbf{x}_k + \mathbf{s})\| < \kappa_{eg}R
\]
and
\[
|f(\mathbf{x}_k + \mathbf{s}) - q(\mathbf{x}_k + \mathbf{s})| < \kappa_{ef}R^2,
\]
for the current iterate $\mathbf{x}_k$ and some constants $\kappa_{eg}, \kappa_{ef} > 0$.

A polynomial interpolation model is fully linear on a set $B$ if the interpolation set $\{ x_1, x_2, \dots, x_m\} = X$ is $\Lambda$-poised on $B$ for a fixed $\Lambda$~\cite[p.~92]{conn_introduction_2009}. A set is $\Lambda$-poised when~\cite[p.42]{conn_introduction_2009}
\[\Lambda \geq \max_{1 \leq i \leq m} \max_{x \in B} |l_i(x)|,\] where the $l_i(x)$ are the Lagrange polynomials, with $\kappa_{ef}, \kappa_{eg}$ depending only on $\Lambda, n, m$, and the Lipschitz constant of $\nabla f$ \cite[p.~92]{conn_introduction_2009}. 

When computing an LU factorization with partial pivoting of the Vandermonde matrix of monomials up to quadratic order for the set $X$, scaled and shifted to lie within the unit ball around the origin, if all pivots $p_i \geq \tau$, then $X$ is $\Lambda$-poised for some $\Lambda(n, \tau)$~\cite[p.~51]{conn_introduction_2009}. Thus, to ensure poisedness, the Conn algorithm builds the set by adding points as long as pivots exceed the fixed threshold $\tau$. CLUSTER-Conn uses this same threshold: it only adds a candidate point $\mathbf{x}_{\text{new}}$ if $p_{m+1} > \tau$ (in line 18). This ensures that, if Algorithm \ref{alg:evaluate_tr} starts with a fully linear model, the returned model is also fully linear. 

\subsection{Additional criticality step}

In addition to constructing a fully linear model when approaching convergence, the criticality step ensures \[\Vert g\Vert \ge \min\{ \epsilon_{\text{c}}, \mu^{-1} R\},\]
which is used to show $R \to 0$~\cite[p.~188]{conn_introduction_2009}. However, because CLUSTER-Conn also adds evaluation points to the model during intermediate (single-cluster) trust-region steps, we need to ensure that this criticality guarantee is not spoiled by more than a constant factor $a$. Thus, we add an additional check at the end of the trust-region step that, at worst, reverts to the original Conn algorithm:
\[\Vert g\Vert < a\min\{ \epsilon, \mu^{-1} R\} \Rightarrow \text{return original point } \mathbf{x}_\text{pot}^{\text{orig}} \text{ and model } (g_0, H_0) \]
for some $a \in (0,1]$. In practice, we set $a=1$ and find that this test is triggered so rarely that it makes little practical difference.

\subsection{Sufficient decrease}

Standard trust-region steps guarantee~\cite{conn_introduction_2009,conn_global_2009, powell_convergence_2012} that
\[ q(\mathbf{x}_\text{k}) - q(\mathbf{x}_\text{k+1}) \geq \kappa_{\text{cd}}[q(\mathbf{x}_\text{k}) - q(\mathbf{x}_\text{k} + \mathbf{s}^C)],\]
where $\mathbf{s}^C$ is the Cauchy step: the minimizer of the quadratic model along the gradient direction and $\kappa_{\text{cd}}>0$ is some constant. 

In CLUSTER, the incremental construction of the step does not preserve this guarantee: in extreme cases, the final point $\mathbf{x}_{\text{new}}$ may even have a higher model value than at $\mathbf{x}_{\text{c}}$ (because the actual objective decreased). Empirically, these steps can accelerate convergence when the model $q$ poorly approximates $f$. 

However, for the convergence guarantee we require the sufficient decrease condition, and so we insert an extra check into the trust-region step of CLUSTER-Conn. Letting $\mathbf{x}^*$ be the minimizer of the model $q$ to be returned, we enforce in CLUSTER-Conn
\[ q(\mathbf{x}_\text{c}) - q(\mathbf{x}_\text{new}) \leq \kappa_3 q(\mathbf{x}_\text{c}) - q(\mathbf{x}^*) \Rightarrow \text{ return } \mathbf{x}^*,\]
with $0 < \kappa_3 \ll 1$.
This is potentially expensive because evaluating the objective at $\mathbf{x}^*$ involves changing all but the last cluster in the loop. By choosing $\kappa_3$ to be small (we choose $\kappa_3 = 10^{-2}$), however, we find that this sufficient-decrease condition is triggered rarely enough to have little impact on performance.

\begin{figure}[t!]
    \centering
    \includegraphics[width=0.75\linewidth]{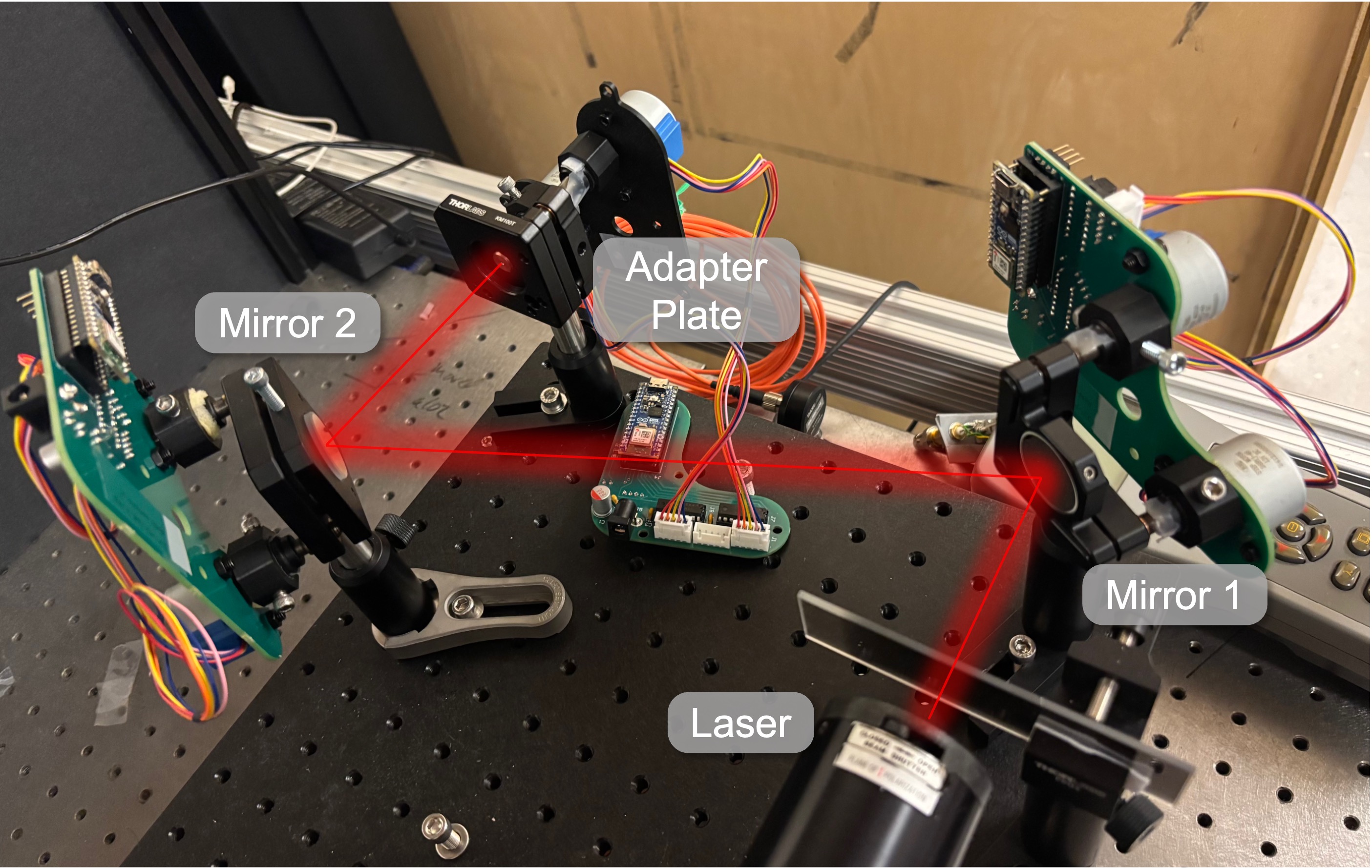}
    \label{fig:fig6}
    \caption{The laboratory experiment used to compare NEWUOA and CLUSTER-Powell.}
\end{figure}

\section{Experimental details}

The laboratory experiment consists of an attenuated $633$~nm Helium-Neon laser beam reflecting from two protected silver mirrors to enter a multimode graded-index optical fiber (model OM1, NA=0.275, diameter = 62.5 $\mu$m) secured in a fiber adapter plate (model SM1FC2). A Thorlabs PM120 power meter measured the coupled power at the fiber output. The configuration is shown in Fig.~\ref{fig:fig6}. For easier comparison between the original NEWUOA algorithm and CLUSTER-Powell, we attached a motor controller to each knob instead of using a robotic arm as in Ref.~\citenum{uddin_ai-driven_2025}. A description of this custom fine-adjustment tool is given in Ref.~\citenum{choi_framework_2026}.

In the Supplementary Information, to demonstrate the applicability of CLUSTER to other experimental tasks, we also tested CLUSTER-Powell on a Michelson interferometer~\cite{michelson_relative_1887}, aiming to align two reflected beams at a specified target pixel on the detector (camera). In this situation, the mirror tip/tilt parameters are decoupled in the objective function, so intermediate evaluations offer little information and CLUSTER updates do not significantly speed up convergence. Nonetheless, we find that CLUSTER-Powell performs as well as or slightly better than NEWUOA.

\end{appendix}

\section*{Acknowledgments}

This work was supported in part by the MIT Generative AI Impact Consortium, and by the U.S. Army Research Office (ARO) through the Institute for Soldier Nanotechnologies (ISN) under award no.~W911NF-23-2-0121. S.L. acknowledges support from an MIT Dean of Science Fellowship. S.P. acknowledges the financial support of the Hertz Fellowship Program and NSF Graduate Research Fellowship Program.

\bibliographystyle{siamplain}
\bibliography{refs.bib}
\end{document}


\maketitle
\footnotetext[1]{Department of Physics, Massachusetts Institute of Technology.}
\footnotetext[2]{Department of Electrical Engineering and Computer Science, Massachusetts Institute of Technology.}
\footnotetext[3]{Nokia Bell Labs.}
\footnotetext[4]{Research Laboratory of Electronics, Massachusetts Institute of Technology.}
\footnotetext[5]{NSF Institute for Artificial Intelligence and Fundamental Interactions (IAIFI).}
\footnotetext[6]{Department of Mathematics, Massachusetts Institute of Technology.}

\begin{figure}
    \centering
    \includegraphics[width=0.75\linewidth]{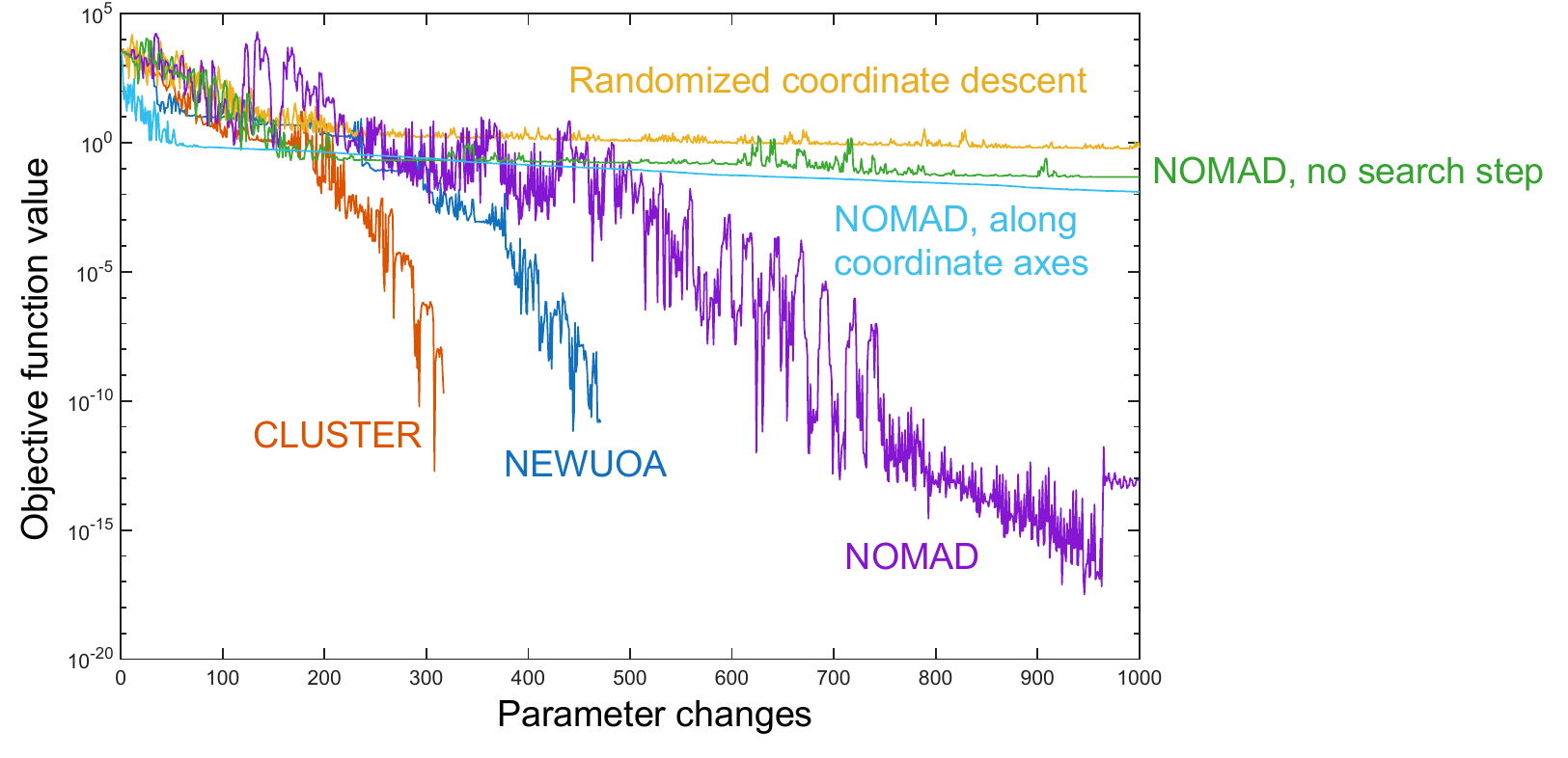}
    \label{fig:coordinate}
    \caption{The convergence history for the helical valley function. }
\end{figure}

As discussed in the main text, we compared CLUSTER-Powell and NEWUOA with several coordinate-search methods. In Fig.~\ref{fig:coordinate}, we show the convergence history for the helical valley function. In addition to the simple randomized coordinate descent algorithm, we tested the NOMAD algorithm~\cite{audet_algorithm_2022}. We tested the standard version of NOMAD, which includes a search step (shown in the figure in purple), as well as a version without the search step (green), and only along the coordinate axes (light blue). 

\begin{figure}[t!]
    \includegraphics[width=\linewidth]{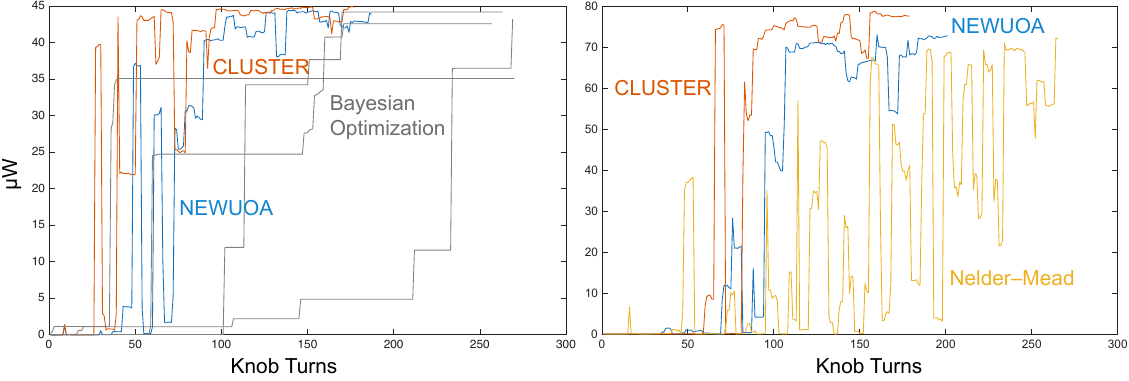}
    \label{fig:bo_nm}
    \caption{Comparison of convergence histories of NEWUOA and CLUSTER-Powell with Bayesian Optimization (BO) (left) and Nelder--Mead (right). For BO, the best-so-far value is plotted, and we show the convergence history of three trials with different sets of randomly-selected initialization points.}
\end{figure}

We also attempted Bayesian optimization and Nelder--Mead on the fiber coupling laboratory experiment. A representative convergence history for each algorithm, compared with NEWUOA and CLUSTER-Powell, is given in Fig.~\ref{fig:bo_nm}.

\begin{figure}[t!]
    \centering
    \includegraphics[width=0.75\linewidth]{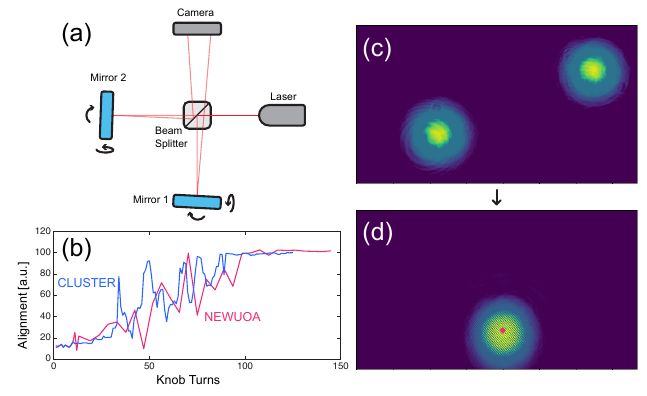}
    \label{fig:interferometer}
    \caption{\emph{Left:} The Michelson interferometer configuration and a representative convergence history. \emph{Right:} The camera view before and after optimization. The target pixel is marked in pink. }
\end{figure}

In addition to the fiber coupling configuration, we also tested a Michelson interferometer experiment~\cite{michelson_relative_1887}. The configuration is shown in Fig.~\ref{fig:interferometer}, as well as a representative convergence history. The parameters are decoupled in this experiment, so the intermediate evaluations of CLUSTER do not offer significant information and we find that NEWUOA and CLUSTER-Powell typically have the same speed of convergence.

\bibliographystyle{siamplain}
\bibliography{refs.bib}